%% file: NortonRoyset_DRM_revision.tex
\documentclass[titlepage,final,11pt]{article}

\usepackage{subcaption}
\usepackage{hyperref}
\usepackage{setspace}
\input mac.tex

\begin{document}


\begin{center}
\begin{large}
{\bf Diametrical Risk Minimization: Theory and Computations}
\smallskip
\end{large}
\vglue 0.7truecm
\begin{tabular}{lcl}
  \begin{large} {\sl Matthew D. Norton 
                                  } \end{large} & \ \ {\phantom{\&}} \ \ &
   \begin{large} {\sl Johannes O. Royset
   				  } \end{large} \\
  \\
  Target Corporation  && Naval Postgraduate School \\
  mdnorto@gmail.com && joroyset@nps.edu
\end{tabular}

\vskip 0.2truecm

\end{center}

\vskip 1.3truecm

\noindent {\bf Abstract}. \quad The theoretical and empirical performance of Empirical Risk Minimization (ERM) often suffers when loss functions are poorly behaved with large Lipschitz moduli and spurious sharp minimizers. We propose and analyze a counterpart to ERM called Diametrical Risk Minimization (DRM), which accounts for worst-case empirical risks within neighborhoods in parameter space. DRM has generalization bounds that are independent of Lipschitz moduli for convex as well as nonconvex problems and it can be implemented using a practical algorithm based on stochastic gradient descent. Numerical results illustrate the ability of DRM to find quality solutions with low generalization error in sharp empirical risk landscapes from benchmark neural network classification problems with corrupted labels.\\

\vskip 0.5truecm

\halign{&\vtop{\parindent=0pt
   \hangindent2.5em\strut#\strut}\cr
{\bf Keywords}:  empirical risk minimization, generalization error,  solution stability.
                         \cr

Correspondence to: J.O. Royset

{\bf Date}:\quad \ \today \cr}

\baselineskip=15pt

\section{Introduction}

\doublespacing
In stochastic optimization, the minimum value of the empirical risk exhibits a downward bias and the corresponding minimizers are therefore often  poor in terms of their true (population) risk. Lipschitz continuity\footnote{We recall that a function $f$ is {\em Lipschitz continuous} on a set $C$ if there is a positive constant $\kappa$, called the {\em Lipschitz modulus}, such that $|f(w)-f(w')| \leq \kappa\|w-w'\|$ for all $w,w'\in C$.} is often brought in as a critical component in attempts to assess the quality of such minimizers, with the Lipschitz moduli of loss functions relative to model parameters (weights) entering in generalization bounds and other results for Empirical Risk Minimization (ERM) problems; see for example \cite{shalev2010learnability,hardt2016train,charles2018stability,bousquet2002stability, bartlett2017spectrally}. In this work, we propose a counterpart to ERM called Diametrical Risk Minimization (DRM) that possesses a generalization bound which is independent of Lipschitz moduli for convex as well as nonconvex loss functions. Preliminary simulations on benchmark Neural Network (NN) classification problems with MNIST and CIFAR-10 datasets support the hypothesis that when problems have large Lipschitz moduli, DRM is able to locate quality solutions with low generalization error, while ERM comparatively struggles.

\begin{figure}
\begin{subfigure}[b]{.5\textwidth}
\includegraphics[scale=.25]{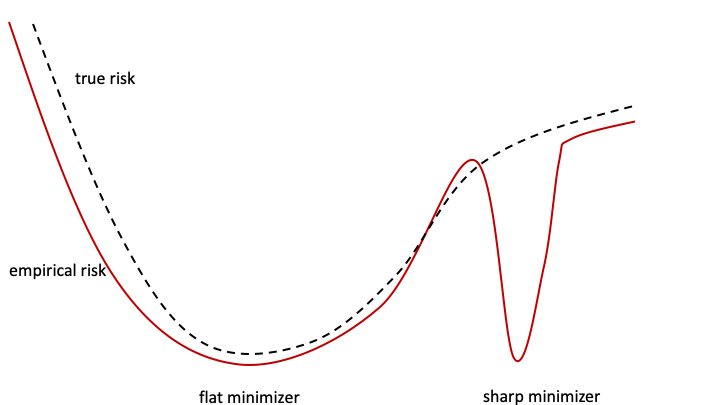}
\caption{Sharp vs flat minimizers.}
\label{fig:sharp_ex_a}
\end{subfigure}
\begin{subfigure}[b]{.5\textwidth}
\includegraphics[scale=.25]{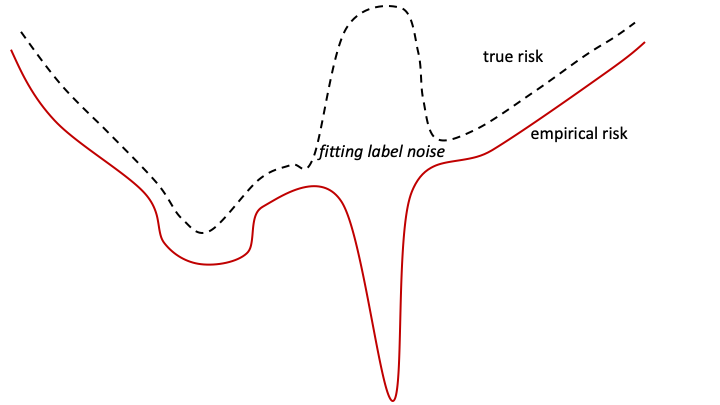}
\caption{Sharp minimizer from label noise in training set.}
\label{fig:sharp_ex_b}
\end{subfigure}
\caption{}
\label{fig:sharp_ex}
\end{figure}

The empirical risk as a function of model parameters in a learning problem has a graph, which we refer to as the {\em empirical risk landscape}; see the solid red line in Figure~\ref{fig:sharp_ex_a}. The process of training the model then amounts to determining parameters that correspond to the bottom of a ``valley'' in this landscape.
A large Lipschitz modulus tends to produce a {\em sharp} empirical risk landscape, where the empirical risk is highly variable with sudden ``drops'' of the kind labelled as a {\em sharp minimizer} in Figure~\ref{fig:sharp_ex_a}. If the Lipschitz modulus is low, at least locally, then the empirical risk landscape is {\em flat} as to the left in Figure~\ref{fig:sharp_ex_a} and the resulting minimizer is {\em flat}.

Instead of the empirical risk, DRM considers the \textit{diametrical risk} at a point in the parameter space, which is given by the worst-case empirical risk in a neighborhood of the point. This provides DRM with a broader view of the empirical risk landscape than ERM and results in improved performance when the landscape is sharp.

Dealing with empirical risk landscapes that have large Lipschitz moduli and sharp minimizers is a challenge that has seen renewed attention in recent years under the heading of sharp vs flat minimizers in landscapes
generated by NN. It is hypothesized that the landscape of NN problems are chaotic \cite{nguyen2017loss, li2018visualizing} and that flat minimizers have better generalization properties compared to sharp ones \cite{keskar2016large,sagun2016eigenvalues,hochreiter1997flat,chaudhari2017entropy}. The potential effects can be seen in Figure~\ref{fig:sharp_ex_a}, which also depicts the true risk in a learning problem (black dashed line).  The spurious dip (right-most valley) of the empirical risk landscape is caused by a large Lipschitz modulus. Since ERM seeks out such dips, the resulting minimizer is poor when assessed using the true risk. We would have preferred that ERM found the left-most valley where any of its obtained minimizers would have had a low true risk. Although the figure is conceptual, it is believed that landscapes of NNs may exhibit similar behavior \cite{keskar2016large}. Thus, it has been a goal of many researchers to either locate the flat minimizers of such problems or to construct loss functions and/or NN architectures which do not have a large number of sharp minimizers \cite{keskar2016large,sagun2016eigenvalues,chaudhari2017entropy,gouk2018regularisation}.

Figure~\ref{fig:sharp_ex_b} illustrates a different kind of sharpness, which can be induced by introducing label noise as explored in \cite{zhang2016understanding,oymak2019generalization}. When a portion of labels are randomly flipped, it has been shown that NNs are capable of fitting the training data perfectly, achieving a zero-training error solution. However, this solution clearly will not generalize and lies near a sharp minimizer which is associated with the fitting of the incorrectly labeled training data. Nevertheless, it has been shown that, even in the presence of label noise, there still exist flat minimizers such as the left-most one in Figure~\ref{fig:sharp_ex_b}. It is immediately clear that DRM may perform well and achieve flat minimizers even in landscapes with spurious sharp minimizers due its broader view of the empirical risk landscape.

Lipschitz moduli frequently play a direct role in statistical learning theory. For example, a paradigm in learning theory is the analysis of algorithmic stability \cite{shalev2010learnability,hardt2016train,charles2018stability,bousquet2002stability}. However, a majority of these learning bounds require some notion of smoothness in terms of either a Lipschitz continuous objective function and/or Lipschitz continuous gradient. The Lipschitz moduli then enter the resulting generalization bounds and influence the (theoretical) stability of the algorithm used to perform ERM.  The reliance of these and other generalization bounds \cite{bartlett2017spectrally} on the magnitude of the Lipschitz moduli, as well as the growing support for the sharp-vs-flat hypothesis, have even given rise to research centered on Lipschitz regularization \cite{oberman2018lipschitz,gouk2018regularisation,qian2018l2} for improving the generalization of NNs. We provide generalization bounds, however, that do not rely upon the Lipschitz moduli: DRM can be applied to sharp empirical risk landscapes with resulting solutions having low generalization error. Eliminating the dependence on the Lipschitz modulus does, however, come with a cost. The provided bounds rely somewhat unfavorably on the dimension of the parameter space. Still, experiments indicate that this is a limitation of the proof approach and not a fundamental limitation of DRM in general. We carry out experiments in a NN setting where the dimension of the parameter space is larger than the number of training samples. Even in this setting, we find that DRM solutions generalize favorably compared to solutions found via ERM.

Although not studied in detail, DRM may also support training in the context of quantification where parameters are represented with lower precision \cite{DongYaoGholamiKeutzerMahoney.19}. There, the explicit robustification in DRM against parameter perturbations could emerge as beneficial.

The downward bias associated with ERM has been known since the early days of stochastic optimization and $M$-estimators. Traditional remedies include a variety of regularization schemes, focused on alteration of the objective function \cite{bousquet2002stability,LiuWangYaoLiYe.18,BertsimasCopenhaver.18} or the optimization procedure itself with, for example, early stopping \cite{hardt2016train,Royset.13,RoysetSzechtman.13}. Another remedy is to replace ERM by the problem of minimizing the distributionally worst-case empirical risk; see for example \cite{WiesemannKuhnSim.14,ZhangXuZhang.16,RoysetWets.16b,DuchiGlynnNamkoong.18,BertsimasGuptaKallus.18} and references therein. Typically, the worst-case is defined in terms of a ball in some metric on a space of probability measures centered on the empirical distribution generated by the available data. Adversarial training \cite{Madry.18} is a closely related approach where the worst-case is computed by perturbing the data directly as for classical robust M-estimators in statistics; see also \cite{liao2018defense,zhang2019theoretically,cohen2019certified,carmon2019unlabeled}
for other efforts to achieve NNs that are robust against data perturbations.
DRM is distinct from these approaches by {\em perturbing the parameter vector} instead of the distribution governing the data or the data itself.

Perturbation of a parameter vector as a means to account for ``implementation error'' of a decision or design specified by the vector is included as a motivation for Robust Optimization \cite{BenTalNemirovski.98}; see \cite{StinstradenHertog.08} for application in the context of meta-models and \cite{MenFreundNguyenSaaseoanePeraire.14} for fabrication problems. The latter reference as well as \cite{Lewis.02,LewisPang.10} lay theoretical and computational foundations for minimizing functions of the form $w\mapsto \nsup_{v\in V} f(w+v)$ that include establishing Lipschitz continuity even if $f$ is rather general. In particular, the minimization of such sup-functions can be achieved by semidefinite programming \cite{Lewis.02,LuoSturmZhang.04} when $f$ is convex and quadratic. Examples of ``robustification'' by considering a worst-case parameter vector are also found in signal processing \cite{Luo.03,PinarArikan.04}.
Concurrent to the present work, \cite{wu2020adversarial} develops an approach based on perturbation of parameters {\em and } data with strong empirical performance. The theoretical results are limited to a PAC-Bayes bound where an assumption on the distribution of parameters allows one to conclude that the approach has a robust generalization bound that involves the expectation of the ``flatness'' of the empirical risk landscape. Though, details of the argument is omitted. A more detailed theoretical analysis is carried in \cite{TsaiHsuYuChen.21} for the specific class of feed-forward neural networks with non-negative monotone activation functions against norm-bounded parameter perturbations. In contrast, we consider nearly arbitrary neural networks and in fact stochastic optimization more generally as well.

The remainder of this paper is organized as follows. Section~\ref{sec:drm} introduces DRM and illustrates the difficulties faced by ERM when loss functions are poorly behaved. In Section~\ref{sec:gen_error}, we provide a theoretical analysis of DRM that includes generalization bounds independent of Lipschitz moduli. In Section~\ref{sec:alg}, we propose a practical algorithm for performing DRM based on stochastic gradient descent (SGD). We then provide an experimental study in Section~\ref{sec:exp}, with a focus on supporting the idea that DRM can find good solutions to problems with sharp empirical risk landscapes. Code for the experiments is available online\footnote{ \url { https://github.com/matthew-norton/Diametrical_Learning } }.

\section{Diametrical Risk Minimization}
\label{sec:drm}
For a loss function $\ell:\reals^n \times\reals^d \to \reals$ and sample $S=\{z_1, \dots, z_m\}\subset \reals^d$, it is well known that the {\em ERM problem}
\[
\nnmin_{w\in \reals^n} R_m(w) = \frac{1}{m}\sum_{i=1}^m \ell(w,z_i)
\]
can be a poor surrogate for the actual problem of minimizing the {\em true risk} $R(w) = \Ex_z[\ell(w,z)]$. Here, $R_m(w)$ is the {\em empirical risk} of parameter vector $w$. Specifically, a (near-)minimizer $w^\star_m$ of $R_m$ tends to have true risk $R(w^\star_m)$ significantly higher than that stipulated by $R_m(w^\star_m)$; cf. the difference between the solid red line illustrating $R_m$ and the dashed black line illustrating $R$ in Figure~\ref{fig:sharp_ex_a}. The effect worsens when the loss $\ell(w,z)$ varies dramatically under changing parameters $w$, which is the case when $\ell(\cdot,z)$ has a large Lipschitz modulus.

In this section, we propose an alternative that we coin {\it Diametrical Risk Minimization (DRM)}. In contrast to common robustification strategies based on perturbing the data set, DRM perturbs the parameters and thereby obtains stability even for poorly behaved loss functions. Instead of minimizing $R_m$ directly as in ERM, any learned parameter vector $w$ is ``diametrically'' modified before the empirical risk is minimized.

\begin{definition}
For a loss function $\ell:\reals^n \times\reals^d \to \reals$ and sample $S=\{z_1, \dots, z_m\}\subset \reals^d$, the {\em diametrical risk} of a parameter vector $w\in\reals^n$ is given as
\[
R_m^\gamma(w) = \sup_{\|v\|\leq \gamma} R_m(w + v) = \sup_{\|v\|\leq \gamma} \frac{1}{m}\sum_{i=1}^m \ell(w+v,z_i),
\]
where $\gamma \in [0,\infty)$ is the {\em diametrical risk radius}.
\end{definition}

We see that the diametrical risk of parameter vector $w\in \reals^n$ is the worst possible empirical risk in a neighborhood of $w$. Any norm $\|\cdot\|$ can be used to define the neighborhood. Trivially, $R^0_m(w) = R_m(w)$, but generally $R^\gamma_m(w) \geq R_m(w)$.

For some set $W\subset\reals^n$ of permissible parameter vectors, the {\it DRM problem} amounts to
\[
\nnmin_{w\in W} R_m^\gamma(w)
\]
and results in a solution $w_m^\gamma$, which might be a global minimizer, a local minimizer, a stationary point, or some other parameter vector with ``low'' diametrical risk.

As we show in Theorem~\ref{tRateofConv}, under mild assumptions,
\begin{equation}\label{eqn:generror_informal}
R(w_m^\gamma) - R_m^\gamma(w_m^\gamma) \leq \beta m^{-1/2}
\end{equation}
with high probability for some constant $\beta$ regardless of the exact nature of $w_m^\gamma$. In particular, $w_m^\gamma$ generalizes even if obtained after aggressive minimization of the diametrical risk; DRM is inherently resistant to overfitting.

Two examples illustrate the challenge faced by ERM when loss functions have large Lipschitz moduli (with respect to parameters $w$). In both examples, we will see that the generalization error for DRM is dramatically smaller than for ERM. For $\kappa\in (1,\infty)$ and $\gamma\in (0,1)$, let
\[
\ell(w,z) = \begin{cases}
\kappa w/\gamma+\kappa          & \mbox{ if } w\in [-\gamma, 0), z = 0\\
-\kappa w/\gamma -\kappa         & \mbox{ if } w\in [-\gamma, 0), z = 1\\
-\kappa w/\gamma+\kappa          & \mbox{ if } w\in [0, \gamma), z = 0\\
\kappa w/\gamma-\kappa         & \mbox{ if } w\in [0, \gamma), z = 1\\
0             & \mbox{ otherwise}.
\end{cases}
\]
If $z$ takes the values 0 and 1, each with probability $\half$, then $R(w) = \Ex_z[\ell(w,z)] = 0$ for all $w\in\reals$. In contrast,
\[
R_m(w) = \begin{cases}
\frac{1}{m} \rho_m (\kappa w/\gamma +\kappa)  & \mbox{ if } w\in [-\gamma, 0)\\
\frac{1}{m} \rho_m (-\kappa w/\gamma+\kappa) & \mbox{ if } w\in [0, \gamma)\\
0             & \mbox{ otherwise}.
\end{cases}
\]
where $\rho_m$ is the number of zeros minus the number of ones in the data $\{z_1, \dots, z_m\}$. Viewing the data as random, we obtain that with probability nearly $\half$, $\rho_m<0$ and thus $w_m^\star = 0$ minimizes $R_m$ for such outcomes of the data and $R_m(w_m^\star) = \rho_m \kappa/m$. Also with probability near $\half$, $\rho_m \geq 0$ and then $w_m^\star = 1$ minimizes $R_m$ and $R_m(w_m^\star) = R(w_m^\star) = 0$. Consequently, $R_m(w_m^\star)$ has a downward bias. Although
\begin{equation}\label{eqn:ex1_error}
R(w_m^\star) - R_m(w_m^\star) \leq \max\{0,-\rho_m\} \kappa/m
\end{equation}
with probability one, the right-hand side includes the constant $\kappa$, which is proportional to the Lipschitz modulus $\kappa/\gamma$ of $\ell(\cdot,z)$. This illustrates the well-known fact that generalization tends to be poor for loss functions with large Lipschitz moduli. However, considering diametrical risk, we have that $R(w_m^\gamma) - R_m^\gamma(w_m^\gamma) \leq 0$ with probability one.

The situation deteriorates further when the loss function is not Lipschitz continuous. Let
\[
\ell(w,z) = \begin{cases}
1/w                              & \mbox{ if } w\in (0,\infty), z = 0\\
-1/w                             & \mbox{ if } w\in (0,\infty), z = 1\\
0                                & \mbox{ otherwise}.
\end{cases}
\]
Again, with $z$ and $\rho_m$ as above, $R(w) = \Ex_z[\ell(w,z)] = 0$ for all $w\in\reals$ and
\[
R_m(w) = \begin{cases}
\frac{1}{m} \rho_m / w   & \mbox{ if } w\in (0,\infty)\\
0                        & \mbox{ otherwise}.
\end{cases}
\]
Then, $\inf_{w\in \reals} R_m(w) = -\infty$ when $\rho_m < 0$, which takes place with probability nearly $\half$. The downward bias is now unbounded. Considering diametrical risk, we have the much more favorable bound $R(w_m^\gamma) - R_m^\gamma(w_m^\gamma) \leq \beta\gamma^{-1} m^{-1/2}$ with high probability for some constant $\beta$.

From these simple examples, it is clear that ERM can lead to arbitrarily slow learning when the loss function is poorly behaved. In the first example above, DRM has a generalization error equal to zero and thus is independent of the Lipschitz modulus. In the second example, DRM reduces the unbounded generalization error encountered by ERM to a quantity proportional to $m^{-1/2}$ as we also see in the following section.

\section{Rates of Convergence}
\label{sec:gen_error}

We begin by formalizing the setting and recall that $f:\reals^n\times \Omega\to \reals$ is a {\em Caratheodory function} relative to a probability space $(\Omega, \cA, P)$ if for all $\omega\in \Omega$, $f(\cdot,\omega)$ is continuous and for all $w\in \reals^n$, $f(w,\cdot)$ is $\cA$-measurable.   In the following, we assume that the data comprises $d$-dimensional random vectors generated by independent sampling according to the distribution $\bbP$ and thus consider the $m$-fold product probability space $(Z^m, \mathcal{Z}^m, \bbP^m)$ constructed from a probability space $(Z, \mathcal{Z}, \bbP)$, with $Z\subset \reals^d$. If $\ell:\reals^n \times Z\to \reals$ is a Caratheodory function relative to $(Z, \mathcal{Z}, \bbP)$, then $R_m$, now viewed as a function on $\reals^n \times Z^m$, is a Caratheodory function relative to $(Z^m, \mathcal{Z}^m, \bbP^m)$; see for example \cite[Prop. 14.44; Ex. 14.29]{VaAn}. Likewise, we view $R_m^\gamma$ as a function on $\reals^n\times Z^m$. It is real-valued by virtue of being the maximum value of the continuous $R_m$ over a compact set\footnote{A set $C\subset \reals^n$ is compact if it is closed and bounded.}. For all $w\in\reals^n$, $R_m^\gamma(w)$ is $\mathcal{Z}^m$-measurable when $R_m$ is a Caratheodory function \cite[Thm. 14.37]{VaAn}. Since $R_m^\gamma$ is continuous (in $w$) for all $(z_1, \dots, z_m)\in \mathcal{Z}^m$ by \cite[Thm. 1.17]{VaAn}, we conclude that $R_m^\gamma$ is a Caratheodory function relative to $(Z^m, \mathcal{Z}^m, \bbP^m)$. In view of \cite[Thm. 14.37; Ex. 14.32]{VaAn}, $\nsup_{w\in W} R_m^\gamma(w)$ is $\mathcal{Z}^m$-measurable as long as $W\subset\reals^n$ is closed. In effect, any concern about measurability in the below statements are put to rest if $\ell$ is a Caratheodory function and $W$ is closed.

We denote by $\Ex$ the expectation with respect to $\bbP$ so that for $w\in\reals^n$, $R(w) = \Ex[\ell(w,z)] = \int \ell(w,z) d\bbP(z)$. When $\ell:\reals^n\times Z$ is a Caratheodory function, we say it is {\em locally sup-integrable} if for all $\bar w\in \reals^n$, there exists $\rho>0$ such that $\int \max\{0, \sup\{\ell(w,z)~|~\|w-\bar w\|\leq \rho\}\} d\bbP(z)<\infty$. It is {\em locally inf-integrable} if ``max-sup'' is replaced by ``min-inf'' in the above statement. The {\em moment-generating function} of a random variable $X$ is $\tau\mapsto E[\exp(\tau X)]$.

We start with a preliminary fact, which follows from Fatou's Lemma.

\begin{proposition}\label{prop:lscusc}
If $\ell:\reals^n\times Z\to \reals$ is a locally inf-integrable Caratheodory function, then $R$ is lower semicontinuous. If locally inf-integrable is replaced by locally sup-integrable, then $R$ is upper semicontinuous.
\end{proposition}
\state Proof. Suppose that $w^\nu\to \bar w$. Since $\ell$ is a inf-integrable Caratheodory function, Fatou's Lemma establishes that $\nliminf \Ex[\ell(w^\nu,z)]\geq \Ex[\nliminf \ell(w^\nu,z)] = \Ex[\ell(\bar w,z)]$. Thus, $R$ is lower semicontinuous. A similar argument confirms the claim about upper semicontinuity; see \cite[Proposition 8.54, 8.55]{primer} for details.\eop

The first main result bounds the amount the diametrical risk can fall below the true risk.

\begin{theorem}{\rm (generalization error in DRM).}\label{thm:generror}
Suppose that $W\subset\reals^n$ is compact, $\ell:\reals^n\times Z\to \reals$ is a locally sup-integrable Caratheodory function, and for all $w\in W$, the moment generating function of $\ell(w,\cdot)-R(w)$ is real-valued in a neighborhood of zero. Then, for any $\epsilon,\gamma>0$ and $m$, there exist $\eta,\beta> 0$ (independent of $m$) such that
\[
\bbP^m\Big(\nsup_{w\in W} \big\{ R(w) - R_m^\gamma(w)  \big\} \leq \epsilon\Big) \geq 1-\eta e^{-\beta m}.
\]
\end{theorem}
\state Proof. By Proposition \ref{prop:lscusc}, $R$ is upper semicontinuous. Let $\{W_k \subset \reals^n, ~k=1, \dots, N\}$ be a finite cover of $W$ consisting of closed balls with radius $\gamma/2$. Without loss of generality, suppose that $W_k \cap W \neq \emptyset$ for all $k = 1, \dots, N$. Let $w^k \in \nargmax_{w\in W_k\cap W} R(w)$, which exists for all $k=1, \dots, N$ because $W_k\cap W$ is nonempty and compact, and $R$ is upper semicontinuous.

For $k=1, \dots, N$, let $\tau\mapsto M_k(\tau)$ be the moment generating function of $R(w^k) - \ell(w^k,\cdot)$ and $I_k(\epsilon) = \nsup_{\tau\in\reals} \{\epsilon \tau - \log M_k(\tau)\}$, which is positive since $M_k$ is real-valued in a neighborhood of zero. Then, by the upper bound in Cramer's Large Deviation Theorem (see for example \cite[Sec. 7.2.8]{ShapiroDenRus.09})
\[
\bbP^m\big( R(w^k) - R_m(w^k) \geq \epsilon \big) \leq e^{-m I_k(\epsilon)}.
\]
Moreover, with $\beta = \min_{k=1, \dots, N} I_k(\epsilon)$,
\[
\bbP^m\Big( \max_{k=1, \dots, N} \big\{R(w^k) - R_m(w^k)\big\} \geq \epsilon \Big) \leq \sum_{k=1}^N  e^{-m I_k(\epsilon)} \leq N e^{-\beta m}.
\]
Consider an event for which $\max_{k=1, \dots, N} \{R(w^k) - R_m(w^k)\} \leq \epsilon$, which takes place with probability at least $1-N e^{-\beta m}$, and let $\bar w\in W$. There exists $\bar k \in \{1, \dots, N\}$ such that $\bar w \in W_{\bar k}$. Since $R_m^\gamma(\bar w) = \nsup_{\|v\|\leq \gamma} R_m(\bar w + v) \geq R_m(w^{\bar k})$,
\[
  R(\bar w) - R_m^\gamma(\bar w) \leq R(\bar w) - R_m(w^{\bar k}) \leq R(\bar w) - R(w^{\bar k}) + \epsilon
  \leq R(w^{\bar k}) - R(w^{\bar k}) + \epsilon  = \epsilon,
\]
where the last inequality follows by the fact that $R(w^{\bar k}) = \nsup_{w\in W_{\bar k}} R(w)$. The conclusion then follows with $\eta = N$ because $\bar w$ is arbitrary.\eop

The theorem furnishes a uniform bound on $R$, which implies in particular that
\[
R(w_m^\gamma) \leq R_m^\gamma(w_m^\gamma) + \epsilon \mbox{ with high probability}
\]
for any parameter vector $w_m^\gamma$ produced by DRM. Thus, there is a strong justification for minimizing $R_m^\gamma$: lower values of the diametrical risk ensure better guarantees on the true risk.
The goal now becomes to develop good methods for producing $w_m^\gamma$ with low $R_m^\gamma(w_m^\gamma)$. The issue of overfitting is largely removed: it is unlikely that a parameter vector $w_m^\gamma$ with a low diametrical risk, i.e.,  low $R_m^\gamma(w_m^\gamma)$, has a high true risk $R(w_m^\gamma)$.

The assumptions in the theorem are generally mild: $\ell(\cdot,z)$ only needs to be continuous and the condition on the  moment generating function is just checked pointwise. The requirement about locally sup-integrable amounts to determine an integrable random variable at every $\bar w$ that upper bounds $\ell$ in a neighborhood of $\bar w$. The constant $\beta$ depends on $\epsilon$, while $\eta$ is a function of $\gamma$ and the diameter of $W$, i.e., $\sup_{w, w'\in W} \|w-w'\|$, in the norm underpinning the diametrical risk. In particular, $\beta$ and $\eta$ are independent of Lipschitz moduli of $\ell$, which may not even be finite.

If the value of the moment generating function of $\ell(w,\cdot)-R(w)$ can be quantified near zero, then we can examine the effect as the error $\epsilon$ vanishes as seen next. We recall that a random variable $X$ is {\em subgaussian} with variance proxy $\sigma^2$ if its mean $E[X]=0$ and its moment generating function satisfies $E[\exp(\tau X)] \leq \exp(\half\sigma^2\tau^2)$ for all $\tau\in\reals$.

\begin{theorem}{\rm (rate of convergence in DRM).}\label{tRateofConv}
Suppose that $W\subset\reals^n$ is compact, $\ell:\reals^n\times Z\to \reals$ is a locally sup-integrable Caratheodory function, and for all $w\in W$, $\ell(w,\cdot)-R(w)$ is subgaussian.
Then, for any $\alpha\in (0,1)$, $\gamma>0$, and $m$, there exists $\beta>0$ (independent of $m$) such that
\[
\bbP^m\Big(\nsup_{w\in W} \big\{ R(w) - R_m^\gamma(w)  \big\} \leq \beta m^{-1/2}\Big) \geq 1-\alpha.
\]
\end{theorem}
\state Proof. By Proposition \ref{prop:lscusc}, $R$ is upper semicontinuous. Let $w^k$ be as in the proof of Theorem \ref{thm:generror}. There exists $\xi\in (0,\infty)$, which may depend on $n$, such that the number of closed balls $N$ of radius $\gamma/2$ required to cover $W$ is no greater than $(\xi/\gamma)^{n}$. Since
$\ell(w^k,\cdot)-R(w^k)$ is subgaussian, say with variance proxy $\sigma_k^2$, we have by Bernstein's
inequality that
\[
\bbP^m\big(R(w^k) - R_m(w^k) > \epsilon\big) \leq \exp\big(-\half m \epsilon^2/\sigma_k^2\big)~~ \mbox{ for all } \epsilon \in [0,\infty).
\]
Let $\sigma = \max_{k=1, \dots, N} \sigma_k$. Thus,
\[
\bbP^m\Big(\max_{k=1, \dots ,N} \big\{R(w^k) - R_m(w^k)\big\} > \epsilon\Big) \leq N\exp\big(-\half m \epsilon^2/\sigma^2\big) \leq \alpha
\]
provided that $\epsilon \geq \beta m^{-1/2}$ and
\[
\beta = \sigma\sqrt{2n \log (\xi/\gamma) - 2\log \alpha}.
\]
Consider the event for which $\max_{k=1, \dots, N} \{R(w^k) - R_m(w^k)\} \leq \epsilon$, which takes place with probability at least $1-\alpha$ for such $\epsilon$. The arguments in  the proof of Theorem \ref{thm:generror} establishes that $\nsup_{w\in W} \{R(w) - R_m^\gamma(w)\} \leq \epsilon$ for this event and the conclusion follows.\eop

The constant $\beta$ in the theorem is given in the proof and depends on the largest variance proxy, denoted by $\sigma^2$, for $\ell(w,\cdot) - R(w)$ at a finite number of different $w$. It also depends on a parameter $\xi$ given by the diameter of $W$. For example, if $R_m^\gamma$ is defined in terms of the sup-norm, then the balls $W_1, \dots, W_N$ in the proof can be constructed according to that norm and the number required is simply\footnote{We note by $\lceil c \rceil$ be lowest integer at least as high as the scalar $c$.} $N =  \lceil \delta/\gamma \rceil^n$, where $\delta = \nsup_{w,\bar w\in W} \|w - \bar w\|_\infty$. This leads to
\[
\beta = \sigma \sqrt{2n\log \lceil \delta/\gamma \rceil-2\log \alpha}.
\]

The constant $\beta$ in the theorem depends unfavorably on $n$. One can attempt to reduce the effect of $n$ by enlarging $\gamma$ as $n$ increases. For example, under the sup-norm one may set
$\gamma = \delta [\exp(\zeta n^{-\nu})]^{-1}$ for positive constants $\zeta$ and $\nu$.  Then, assuming that $\delta/\gamma$ is an integer (which can be achieved by enlarging $\delta$ when needed),
\[
\beta = \sigma\sqrt{2n\log \lceil\delta/\gamma\rceil-2\log \alpha} = \sigma\sqrt{2\zeta n^{1-\nu}-2\log \alpha}.
\]
For example, if $\nu = 1$, then $\beta$ becomes independent of $n$ at the expense of having to grow $\gamma$ rather quickly as $n$ increase. A compromise could be $\nu = \half$, in which case $\beta$ grows only as $n^{1/4}$ and $\gamma$ grows somewhat slower too. Still, in the limit as $n\to \infty$, $R_m^\gamma$ involves maximization over all of $W$, which of course leads to an upper bound.

For fixed $n$, we may also ask what is the right value of $\gamma$. Since a large value might imply additional computational burden and also lead to overly conservative upper bounds, it would be ideal to select it as small as possible without ruining the rate (significantly). One possibility could be to set $\gamma$ proportional to $m^{-1}$ because then the rate deteriorate only with a logarithmic factor from $m^{-1/2}$ to $\sqrt{m^{-1} \log m}$.

It is clear from the proof of the theorem that the assumptions about independent sampling and subgaussian random variables can be relaxed. We only needed that the error in $R_m(w)$ relative to $R(w)$ can be bounded for a finite number of $w$; see \cite{BoucheronLugosiMassart.16,OliveiraThompson.17} for possible extensions.

It may be of interest to determine the error of an obtained parameter vector $w_m^\gamma$ relative to the set of actually optimal parameters $\nargmin_{w\in W} R(w)$.  Theorem \ref{tRateofConv} yields immediately that for any $\delta\in\reals$,
\[
\big\{w\in W~|~R_m^\gamma(w) \leq \delta\big\} \subset \big\{w\in W~|~R(w) \leq \delta + \beta m^{-1/2}\big\} \mbox{ with probability at least } 1-\alpha.
\]

We now examine the harder question of confidence regions for ``good'' parameter vectors. For two sets $A,B\subset\reals^n$, we denote the {\em excess} of $A$ over $B$ by
\[
\exs(A;B) = \sup_{w\in A} \inf_{\bar w\in B} \|w - \bar w\| \mbox{ for nonempty } A, B;
\]
with the convention that $\exs(A;B) = \infty$ if $A\neq\emptyset$ and $B=\emptyset$; $\exs(A;B) = 0$ otherwise. Below, the sets of interest are lower level-sets of $R_m$, possible with a random level. Arguing by means of Prop. 14.33, Thm. 14.37, and Ex. 14.32 in \cite{VaAn}, we see that the excess involving such sets is measurable.

\begin{theorem}{\rm (confidence region in DRM).}\label{tConfReg}
Suppose that $W\subset\reals^n$ is compact, $\ell:\reals^n\times Z\to \reals$ is a locally inf-integrable Caratheodory function, and for all $w\in W$, $\ell(w,\cdot)-R(w)$ is subgaussian. Then, for any $\alpha\in (0,1)$, $\gamma>0$, $\delta\in\reals$, and $m$, there exists $\beta>0$ (independent of $m$) such that
\begin{align*}
\bbP^m\bigg( & \exs\Big( \big\{w\in W~\big|~ R(w) \leq \delta \big\}; ~ \big\{w\in W~\big|~R_m(w) \leq \delta + \beta m^{-1/2}\big\} \Big) \leq \gamma,\\
            & \exs\Big( \nargmin_{w\in W} R(w);  ~\big\{w\in W~\big|~R_m(w) \leq \inf_{\bar w\in W} R_m^\gamma(\bar w) + 2\beta m^{-1/2}\big\} \Big) \leq \gamma\bigg) \geq 1-\alpha.
\end{align*}
\end{theorem}
\state Proof. By Proposition \ref{prop:lscusc}, $R$ is lower semicontinuous. From the compactness of $W\subset \reals^n$, we obtain $\xi\in (0,\infty)$, which may depend on $n$, such that the number of closed balls $N$ with radius $\gamma/2$ required to cover $W$ is no greater than $(\xi/\gamma)^{n}$.
Suppose that $\{W_k \subset \reals^n, ~k=1, \dots, N\}$ is a collection of such balls with $W_k \cap W \neq \emptyset$ for all $k = 1, \dots, N$. Let $w^k \in \nargmin_{w\in W_k\cap W} R(w)$, which exists for all $k=1, \dots, N$ because $W_k\cap W$ is nonempty and compact, and $R$ is lower semicontinuous.

Since $\ell(w^k,\cdot)-R(w^k)$ is subgaussian, say with variance proxy $\sigma_k^2$, we have by Bernstein's
inequality that
\[
\bbP^m\Big(\big|R_m(w^k) - R(w^k)\big| > \epsilon\Big) \leq 2\exp(-\half m \epsilon^2/\sigma_k^2) ~~ \mbox{ for all } \epsilon \in [0, \infty).
\]
Let $\sigma = \max_{k=1, \dots, N} \sigma_k$. Thus,
\[
\bbP^m\Big(\max_{k=1, \dots, N} \big|R_m(w^k) - R(w^k)\big| > \epsilon\Big) \leq 2N\exp(-\half m \epsilon^2/\sigma^2) \leq \alpha
\]
provided that $\epsilon \geq \beta m^{-1/2}$ and
\[
\beta = \sigma\sqrt{2n \log (\xi/\gamma) - 2\log (\alpha/2)}.
\]
Consider the event for which $\max_{k=1, \dots, N} |R_m(w^k) - R(w^k)| \leq \epsilon$. Let $\bar w\in W$ satisfy $R(\bar w)\leq \delta$. Then, there exists $\bar k$ such that $\bar w\in W_{\bar k}$ and because $R(w_{\bar k}) = \inf_{w\in W_{\bar k} \cap W} R(w)$ we obtain that
\[
R_m(w_{\bar k}) \leq R(w_{\bar k}) + \epsilon \leq  R(\bar w) + \epsilon \leq \delta + \epsilon.
\]
Since $\|\bar w - w_{\bar k}\|\leq \gamma$, we conclude that
\[
\exs\big( \{w\in W~|~ R(w) \leq \delta \}; ~ \{w\in W~|~R_m(w) \leq \delta + \epsilon\} \big) \leq \gamma.
\]

We next turn to the result for $\nargmin_{w\in W} R(w)$ and let $w^\star$ be a point in that set. Then, there is $k^\star$ such that $w^\star \in W_{k^\star}$ and
\[
R_m(w^{k^\star}) \leq R(w^{k^\star}) + \epsilon \leq R(w^\star) + \epsilon = \ninf_{w\in W} R(w) + \epsilon.
\]
Moveover, let $\bar w \in \nargmin_{w\in W} R_m^\gamma(w)$. Then, there is $\bar k$ such that $\bar w\in W_{\bar k}$ and
\[
\ninf_{w\in W} R_m^\gamma(w) = \nsup_{\|v\|\leq \gamma} R_m(\bar w+v) \geq R_m(w^{\bar k}) \geq R(w^{\bar k}) - \epsilon \geq \ninf_{w\in W} R(w) - \epsilon.
\]
Combining these inequalities, we obtain that
\[
R_m(w^{k^\star}) \leq \ninf_{w\in W} R_m^\gamma(w) + 2\epsilon.
\]
Since $\|w^\star - w^{k^\star}\|\leq \gamma$, we have established that
\[
\exs\big( \nargmin_{w\in W} R(w);  ~\{w\in W~|~R_m(w) \leq \inf_{\bar w\in W} R_m^\gamma(\bar w) + 2\epsilon\} \big) \leq \gamma
\]
and the conclusion follows.\eop

Since the constant $\beta$ is nearly of the same form as in Theorem \ref{tRateofConv}, the discussion following that theorem carries over to the present context. In particular, we note that a lower-level set of $R_m$, enlarged with $\gamma$, covers $\nargmin_{w\in W} R(w)$ with high probability.

\section{Algorithms for Diametrical Risk Minimization}
\label{sec:alg}
Although there are some computational challenges associated with DRM, most of the existing optimization procedures for ERM can be adapted. Significantly, if the empirical risk function $R_m$ is convex, then the diametrical risk function $R_m^\gamma$ is also convex. Moreover, regardless of convexity,
\[
\frac{1}{m} \sum_{i=1}^m \nabla_w \ell(\bar w+\bar v,z_i)
\]
is a subgradient (in the general sense, cf. \cite[Ch. 8]{VaAn}) of $R_m^\gamma$ at $\bar w$ under weak assumptions,
where $\bar v \in \nargmax_{\|v\|\leq \gamma} R_m^\gamma(\bar w + v)$; see for example \cite[Cor. 10.9]{VaAn}. This implies that standard (sub)gradient methods apply provided that $\bar v$ can be computed. Since $\gamma$ might very well be small, this could be within reach, at least approximately, by carrying one iteration of gradient {\it ascent}. However, this could become costly as computation of such subgradients need to access all data points. This challenge is similar to the one faced by adversarial training \cite{Madry.18}, but there the gradient ascent is carried out relative to the data; see also \cite{ZhengChenRen.18,GongRenYeLiu.20,WongRiceKolter.20}.

We utilize a less costly approach based on the application of SGD to an outer approximation formed via sampling. In short, we approximate the inner maximization by maximizing over a finite set of random points inside the $\gamma$-neighborhood at the current solution $w^t$ of each iteration. We find this approach to be effective, even when working with problems involving NNs where the dimension of $w$ is in the millions. In these applications an outer approximation of $R_m^\gamma$ using as little as 10-20 samples from $\{ v \; | \; \|v\|= \gamma \}$ suffices to achieve improvement over ERM.

We observe also that DRM is related to but distinct from ERM with early termination. If from a minimizer $w_m^\star$ of $R_m^\gamma$ the process of maximizing $R_m(w_m^\star + v)$ subject $\|v\|\leq \gamma$ follows the trajectory along which the algorithm approached $w_m^\star$ in the first place, then DRM would be equivalent to ERM that terminates a distance $\gamma$ from a minimizer. However, this equivalence will only take place when $w_m^\star$ is approached along such direction. It appears that this will occur only occasionally.

%
%
%
%
%
%
%
%
%

\subsection{Gradient Based Algorithm}

We propose two variations of an SGD-based algorithm for DRM which we denote by Simple-SGD-DRM and SGD-DRM. 
 We start with a simple version of the main algorithm that is easier to follow and then introduce the full algorithm with minor alterations aimed toward improving efficiency. In the following, let $\prj_W(w)$ denote the projection of $w$ on $W$ and let $R_{B_t }(w) = \frac{1}{|B_t|}\sum_{z \in B_t} \ell(w,z)$ denote the empirical risk over a batch $B_t \subset S$.\\

\state Algorithm 1: Simple-SGD-DRM.

\state Step 0. Initialize $w^0\in W$,  $r \in \nats$, $t = 0$. Initialize sequence of batches $B_t \subset S$ and learning \\ \indent \indent \; \;\;rates $\lambda_t >0$ for $t=1,\dots,T$.

\state Step 1. Sample $r$ random perturbations (directions): $U=\{u_1, \dots, u_{r}  \; | \; \|u\| = \gamma \}$

\state Step 2. Select $u^\star\in \nargmax_{u \in U}  \frac{1}{|B_t|}\sum_{z \in B_t } \ell(w^t + u, z)$

\state Step 3. Compute $w^{t+1} = \prj_W\big(w^t - \lambda_{t} \nabla_w R_{B_t }(w^t +u^\star) \big)$.

\state Step 4. If $t=T$, stop. Else, $t \leftarrow t+1$ and return to Step 1. \\

The Simple-SGD-DRM algorithm, at each iteration $t$, performs an SGD update towards minimizing the approximating objective function
 \[
w\mapsto \max_{u \in U } \frac{1}{|B_t|}\sum_{z \in B_t }  \ell (w+u , z_i).
 \]
The algorithm does so by first forming a set of $r$ random directions (vectors) $U = \{u_1,\dots,u_{r} \}$ with norm equal to $\gamma$. Then, it determines the more critical $u\in U$, i.e., $u^\star \in \nargmax_{u \in U}  \frac{1}{|B_t|}\sum_{z \in B_t }$ $\ell (w^t +u , z)$. A subgradient of the approximating objective function is then $\nabla_w R_{B_t }(w^t +u^\star )$.

This algorithm, however, does have drawbacks. First, sampling $r$ new vectors in Step 1 at every iteration can be computationally expensive. It may be enough, as we will see in the experiments, to only perform Step 1 intermittently.\footnote{In experiments, we sample only every 5th iteration.} Second, it may be beneficial to reuse one or more of the sampled vectors from Step 1 in future iterations, particularly if we decide to perform sampling only intermittently.

The following algorithm, which we simply call SGD-DRM, includes these options explicitly. One will notice that it can be made equivalent to simple-SGD-DRM with particular choices, and is thus an extension with more options to save computation by limiting sampling. \\

\state Algorithm 2: SGD-DRM.

\state Step 0. Initialize $V_{-1} = \{ \empty \}$, $w^0\in W$, $q\in \nats$,  $r\in \nats$, $p\in[0,1]$, $t = 0$. Initialize sequence of \\ \indent \indent \; \;\;batches $B_t \subset S$ and learning rates $\lambda_t >0$ for $t=1,\dots,T$.

\state Step 1. Sample $r$ random perturbations (directions): $U=\{u_1, \dots, u_r  \; | \; \|u\| = \gamma \}$

\state Step 2. Select $u^\star\in \nargmax_{u \in U}  \frac{1}{|B_t|}\sum_{z \in B_t } \ell(w^t + u, z)$

\state Step 3. Let $V_t = V_{t-1} \cup \{u^\star\} $. If $|V_t| > q$, remove oldest element from $V_t$

\state Step 4. Select $v^\star \in \nargmax_{v \in V_t}  \frac{1}{|B_t|}\sum_{z \in B_t } \ell(w^t + v, z)$

\state Step 5. Compute $w^{t+1} = \prj_W\big(w^t - \lambda_{t} \nabla_w R_{B_t }(w^t +v^\star) \big)$.

\state Step 6. If $t=T$, stop. Else, $t \leftarrow t+1$ and with probability $p$, return to Step 1; with probability \\ \indent \indent \; \;\;$1-p$ return to Step 4 with $V_t = V_{t-1}$. \\

The primary difference between this and the former algorithm is within Steps 3-4 and Step 6. Step 6 allows one to skip the expensive sampling in Step 1 at some frequency represented by $p$. The new set $V_t$ is introduced in Steps 3-4 to allow the reuse of one or more vectors $u$ from previous iterations. As the algorithm progresses, the set $V_t$ acts like a queue with maximal size $q$.\footnote{In experiments, we set $q=1$.} Every time the sampling of Step 1 is not skipped, $V_t$ will be equal to the set $V_{t-1}$ with its oldest element replaced by $u^\star$. For iterations $t\leq q$, the oldest element need not be removed since the queue has not reached its maximum length of $q$. 

\subsection{Implementation for Neural Networks}
\label{sec:drm_alg_nn}
In the following experiments, we consider NN classifiers. Because of the structure of NNs we implement the perturbations $w+v$ and $w+u$ by considering the groupings of parameters that correspond to the structure of each NN layer. For example, a two layer NN might have parameter matrices $\{W_1,W_2\}$, each belonging to the separate network layers. Because of this, we select perturbations in Step 1 such that the {\em Frobenius norm}\footnote{The Frobenius norm of a matrix $A$ is defined as $(\sum_{i,j} a_{ij}^2)^{1/2}$, where $a_{ij}$ are the elements of $A$.} of each layer-wise perturbation matrix is equal to $\gamma$. So, for the network with separable parameters $\{W_1,W_2\}$, a single sample from Step 1 would look like $\{ U_1 , U_2 \}$ with $\|U_1\|_F = \gamma, \|U_2\|_F = \gamma$. Additionally, in the implementation, we first sample each component from a standard normal distribution, then normalize the resulting vector (or matrix) to have norm equal to $\gamma$.

Additionally, in the experiments we implement the coin flip (based on $p$) from Step 6 deterministically. We only perform Step 2 and 3 at every 5th step. Otherwise, we let $V_t = V_{t-1}$. This allows us to save computation time, particularly since the sampling of $U$ can be expensive for NNs with millions of parameters.

\section{Experiments}
\label{sec:exp}
We aim to illustrate that DRM is resistant to overfitting and that its solutions have different local characteristics compared to those from ERM. In particular, we hypothesize that minimizers found by DRM lie in flat portions of the empirical risk landscape. To illustrate these aspects, we focus on the problem of classification with NNs in the presence of label noise, which is suitable because it is prone to overfitting as ERM settles into sharp minimizers; see \cite{oymak2019generalization,zhang2016understanding}. These works show that NNs have the unique ability to perfectly fit (with zero training error) both a data set with label noise as well as the same data with correct labels, with large generalization error in the former but small generalization error in the latter.

In these experiments, we train NNs on a subset of the MNIST and CIFAR-10 datasets with large amounts of label noise (50\% of training labels flipped to an incorrect class). Using standard SGD for ERM (labelled SGD-ERM), the NNs indeed settle into solutions with high generalization error. On the other hand, SGD-DRM is remarkably resistant to overfitting, finding solutions with dramatically lower generalization error. Additionally, we find empirical evidence that the SGD-DRM solution lies in a flat portion of the empirical risk landscape compared to the SGD-ERM solution which appears in a sharper portion.

Code associated with these experiments and a PyTorch\footnote{ \url { https://pytorch.org/ } } based implementation of SGD-DRM is available online\footnote{ \url { https://github.com/matthew-norton/Diametrical_Learning } }. The computational resources is a single Tesla V100 GPU, 16 core Xeon processor and 64 GB memory. The operating system is Ubuntu. Generally, training time for DRM is approximately 3-5 times longer than those for ERM with the same architecture and data. We conjecture that the training times can be reduced by adapting the ideas from fast adversarial training; see for example \cite{ZhengChenRen.18,GongRenYeLiu.20,WongRiceKolter.20}.

 \begin{figure}[t]
 \begin{subfigure}[b]{.5\linewidth}
  \includegraphics[scale=.4]{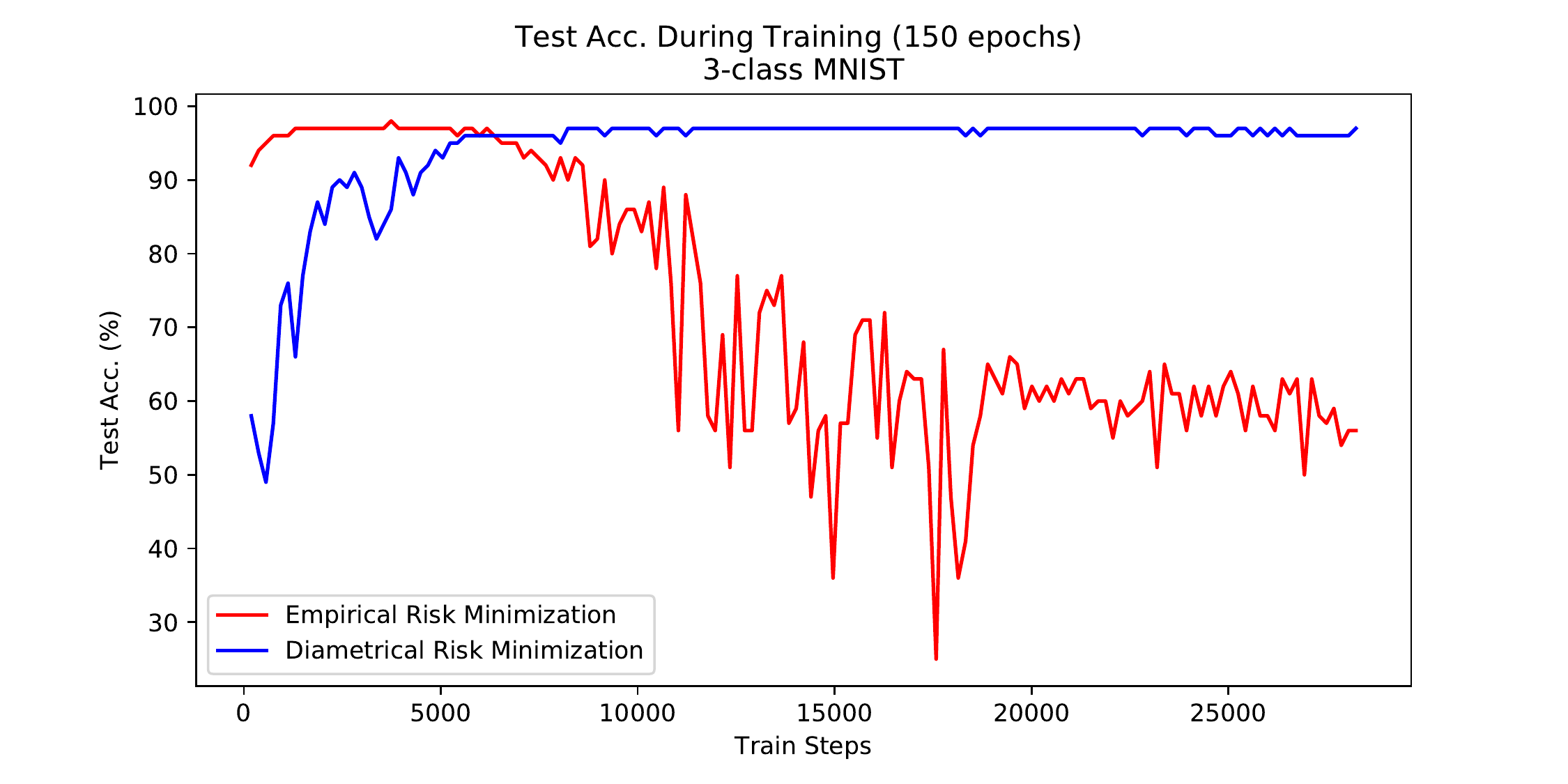}
         \caption{}
    \label{fig:mnist_150_acc}
    \end{subfigure}
 \begin{subfigure}[b]{.5\linewidth}
  \includegraphics[scale=.4]{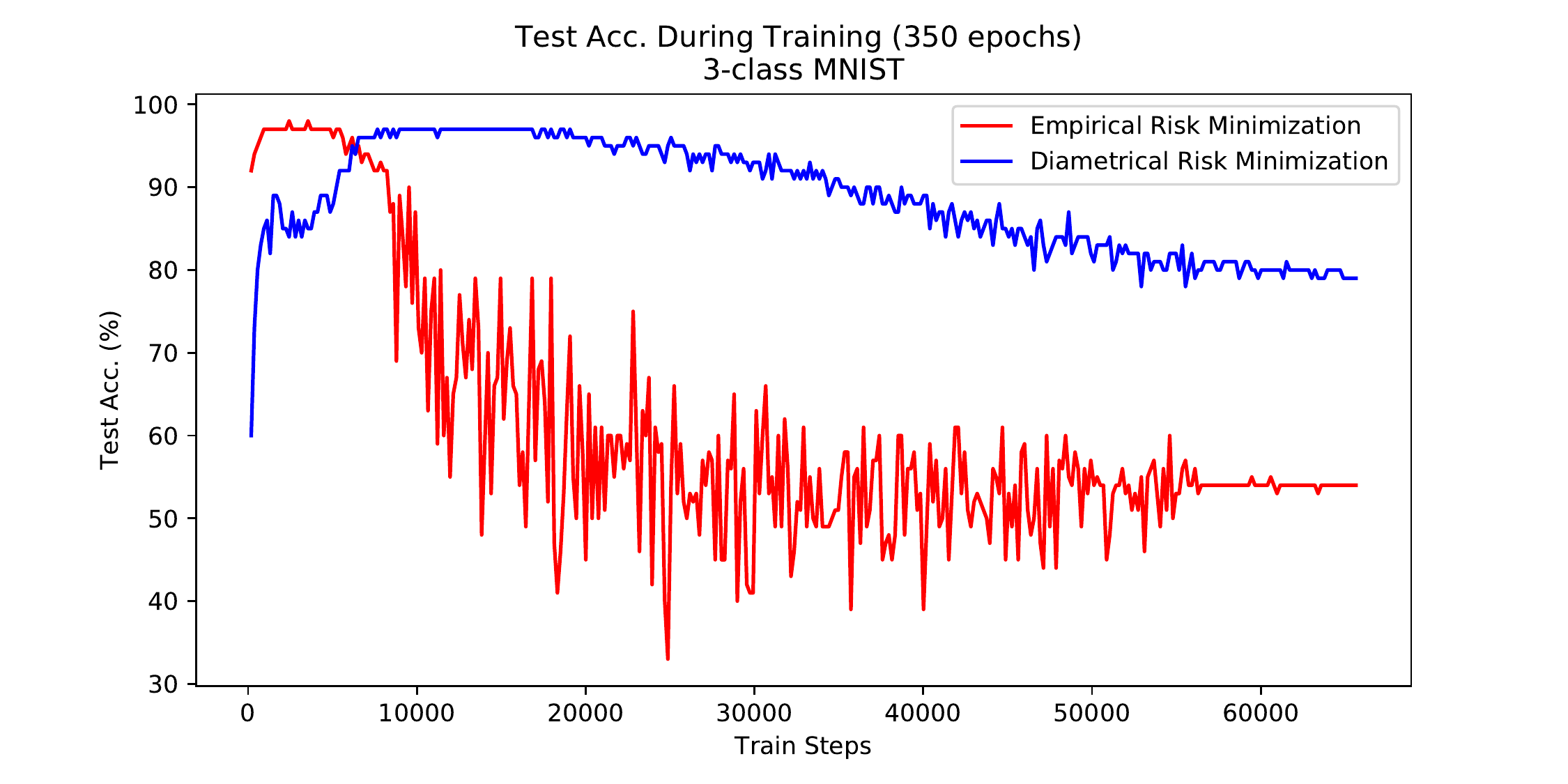}
  \caption{}
  \label{fig:mnist_350_acc}
      \end{subfigure}

 \begin{subfigure}[b]{.5\linewidth}
  \includegraphics[scale=.4]{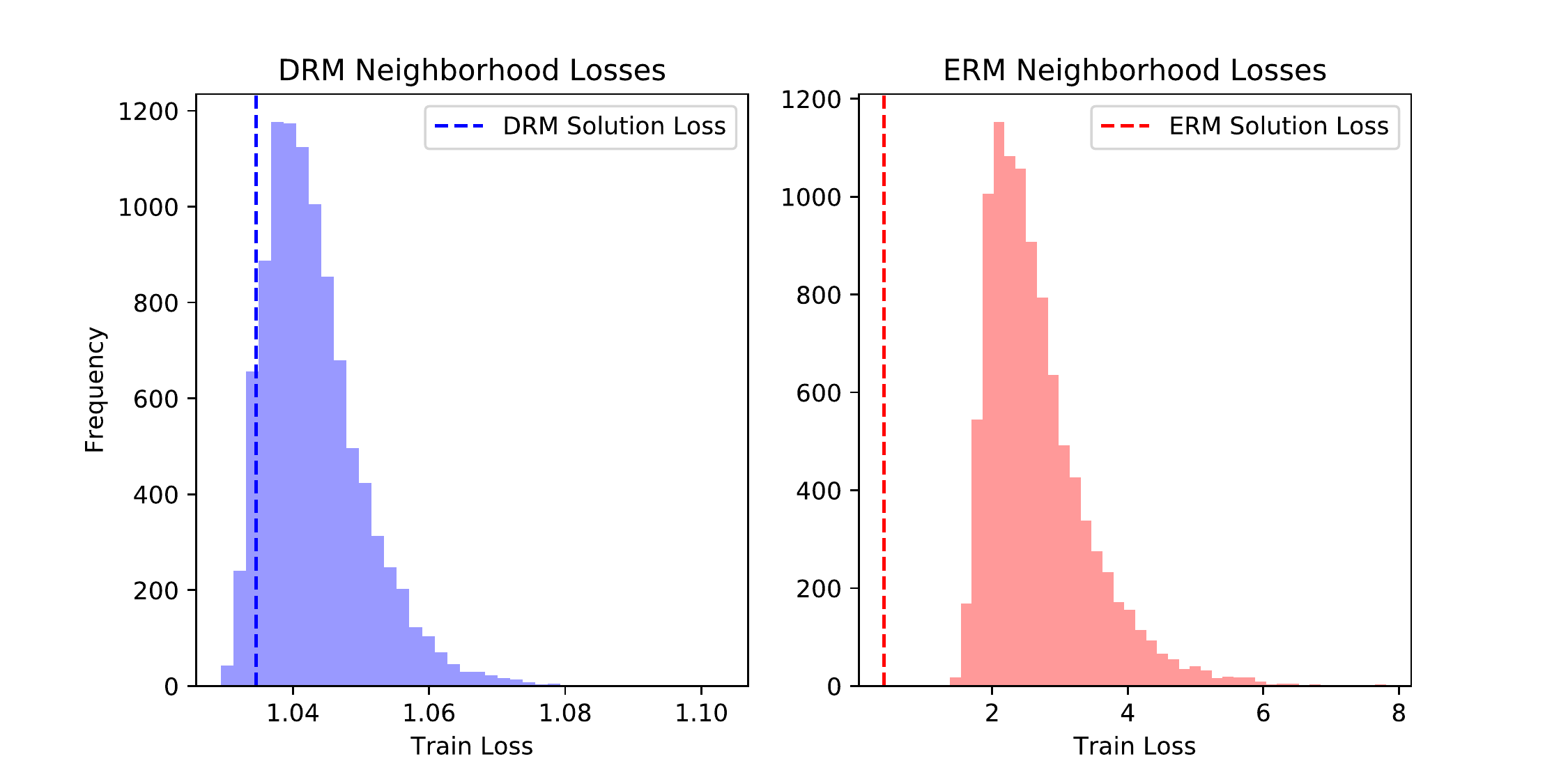}
      \caption{}
           \label{fig:mnist_150_lipschitz}
    \end{subfigure}
 \begin{subfigure}[b]{.5\linewidth}
  \includegraphics[scale=.4]{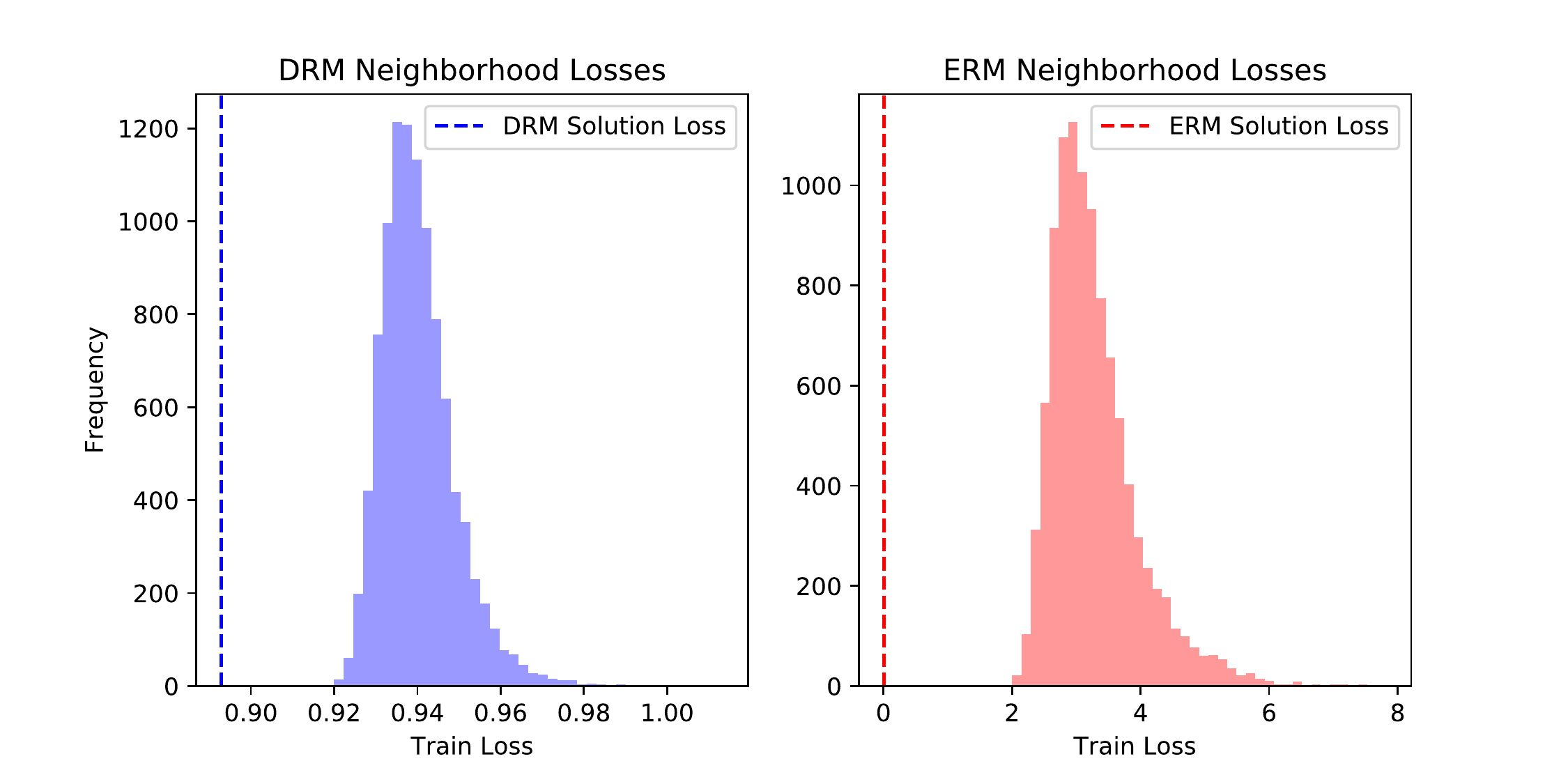}
    \caption{}
      \label{fig:mnist_350_lipschitz}
      \end{subfigure}
            \label{fig:mnist}
\caption{MNIST experiments: (a) Test accuracy when trained for 150 epochs. (b) Test accuracy when trained for 350 epochs. (c) Distribution of empirical risk for 10,000 random points in $\gamma$-neighborhood of final 150 epoch DRM and ERM solution. (d) Distribution of empirical risk for 10,000 random points in $\gamma$-neighborhood of final 350 epoch DRM and ERM solution. }
\end{figure}
\subsection{MNIST}

In the first set of experiments, we train a fully connected, 3 layer NN with hidden units per layer being (320, 320, 200) with ReLU nonlinearities and an additional fully connected output layer feeding into a 3-class softmax negative log-likelihood objective function. For the MNIST dataset, we use only the handwritten digits zero, one, and two so that SGD-ERM achieves nearly zero training error. We flip 50\% of the training labels to an incorrect class. For both SGD-ERM and SGD-DRM, we use standard SGD updates with batch size 100 and learning rate .01 until the last 50 epochs when decreased to .001. For SGD-DRM, the Step 1 perturbations are treated on a layer-wise basis (see Section~\ref{sec:drm_alg_nn}) with $\gamma = 10$. We also use $r=20$ and $q=1$ for the size of $U$ and $V$, respectively. We implement Step 6 deterministically with Step 1 being performed every 5th iteration (batch).

We train twice. First, we train the network for 150 epochs total with learning rate .01 for the first 100 and .001 for the final 50. Test accuracy can be seen in Figure~\ref{fig:mnist_150_acc}. We first note the behavior of SGD-ERM. It begins by finding a good solution that generalizes, but then continues minimizing the empirical risk and settles into a solution that fits the training data with incorrect labels and thus suffers from a sharp decline in test accuracy (a sharp increase in generalization error). In contrast, SGD-DRM resists overfitting. Once it finds a good solution that generalizes, it is able to stay there, resisting the fall into a poor solution. As mentioned earlier, it is hypothesized that solutions in flat portions of the empirical risk landscape generalize better than those in sharp portions. Figure~\ref{fig:mnist_150_lipschitz} plots the distribution of empirical risk for 10,000 random points in the $\gamma$-neighborhood of the final 150-epoch DRM and ERM solutions and illustrates that, indeed, the DRM solution is in a much flatter portion of the empirical risk landscape than the ERM solution. The dashed line represents the value of the empirical risk at the found solution and the rest of the plot is the distribution of empirical risk at points surrounding the solution\footnote{The points are sampled in the same way as for $U$ in Step 1 of SGD-DRM, with neighborhood points sampled on a layer-wise basis as $\{w^* + u \; | \;  \|u\|= \gamma\}$. Additionally, we use the same set of points $u$ for approximating the neighborhood of ERM and DRM solutions.}. Figure~\ref{fig:mnist_150_lipschitz} also illustrates that, while the empirical risk of the SGD-ERM solution is lower (the red dashed line) than that of SGD-DRM it has much larger diametrical risk.

To allow SGD-ERM enough time to achieve nearly zero training error, we also train for 350 epochs with learning rate .01 for the first 300 and .001 for the final 50. Figure~\ref{fig:mnist_350_acc} depicts test accuracy. Again, we see the same behavior for SGD-ERM, as it chaotically falls into a poor solution that does not generalize. SGD-DRM remains resistant to overfitting. While it does experience some degradation in test accuracy, it still ends at a much better solution and its path is not nearly as chaotic; see the smooth vs choppy lines in Figure~\ref{fig:mnist_350_acc}. We also see, again, that the SGD-DRM solution lies in a flatter portion of the empirical risk landscape than that from SGD-ERM. Figure~\ref{fig:mnist_350_lipschitz} illustrates that, while the empirical risk of the SGD-ERM solution is lower (the dashed line at zero) it has (approximately) much larger diametrical risk equal to around 6. The SGD-DRM solution, on the other hand, has higher empirical risk ($\approx$.89) but much smaller diametrical risk ($\approx$.98).

 \begin{figure}[t]
 \begin{subfigure}[b]{.5\linewidth}
  \includegraphics[scale=.4]{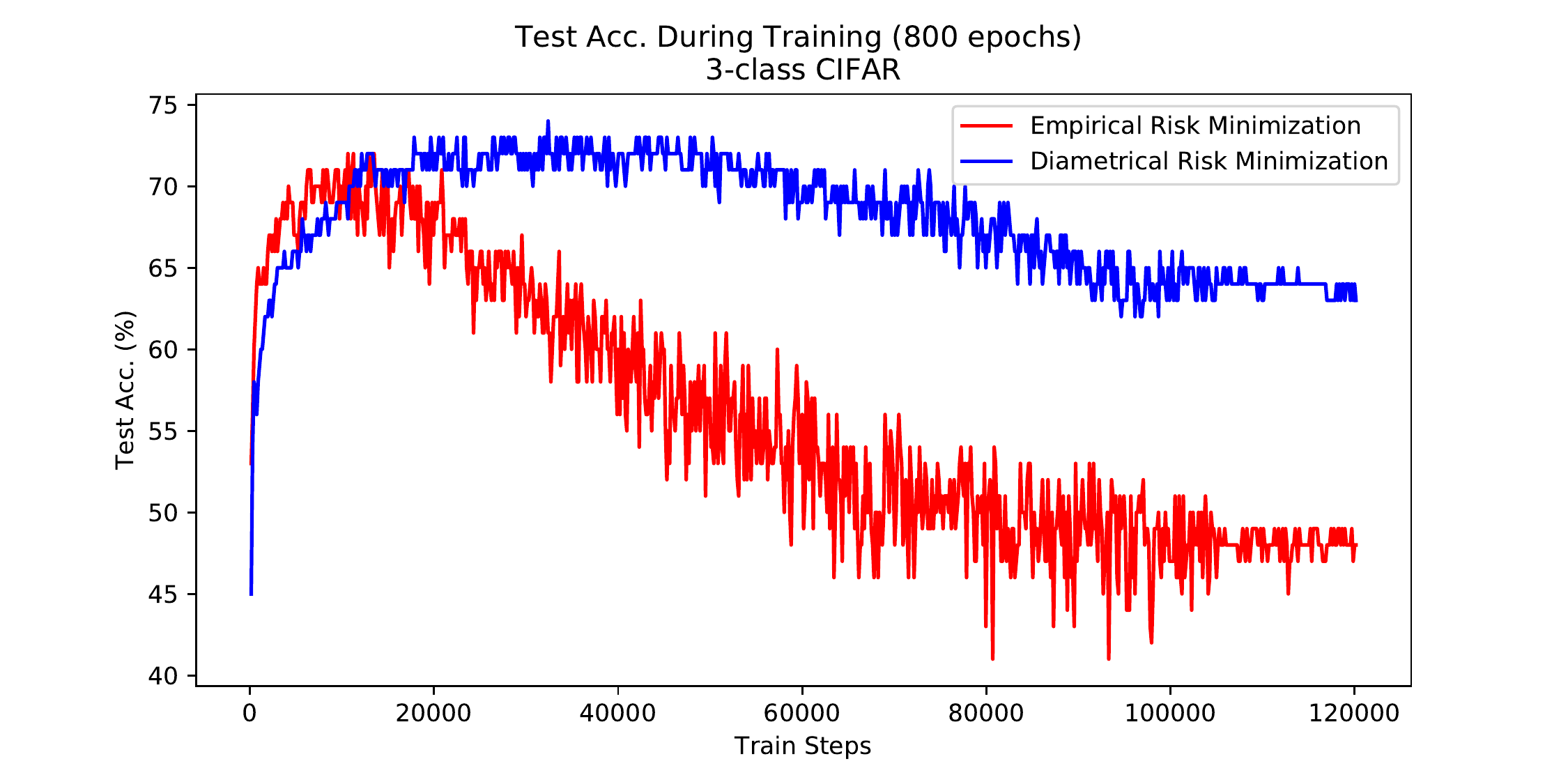}
  \caption{Fully Connected}
    \label{fig:cifar_simple_acc}
    \end{subfigure}
 \begin{subfigure}[b]{.5\linewidth}
  \includegraphics[scale=.4]{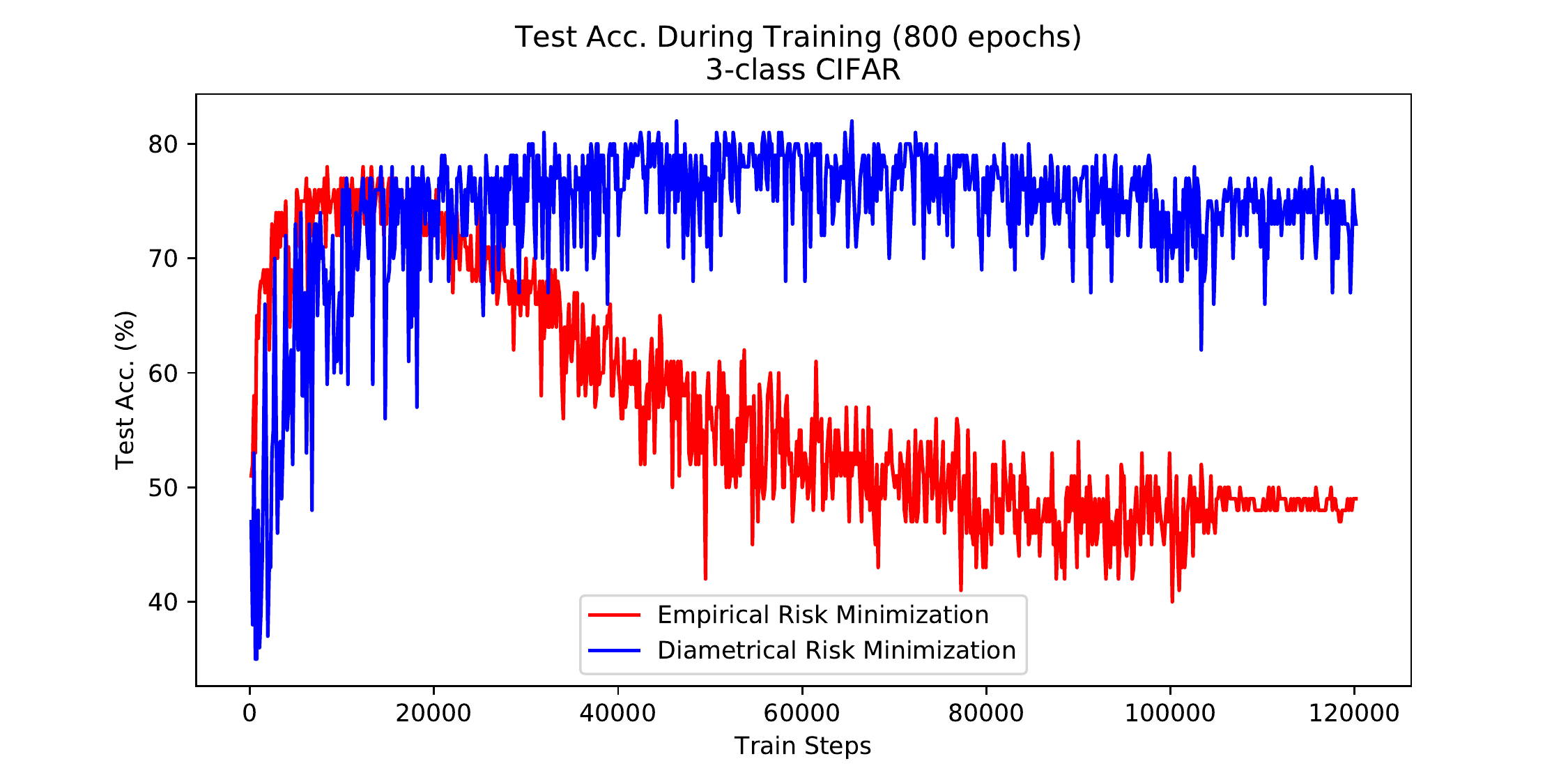}
  \caption{ResNet20}
  \label{fig:cifar_resnet_acc}
      \end{subfigure}

 \begin{subfigure}[b]{.5\linewidth}
  \includegraphics[scale=.4]{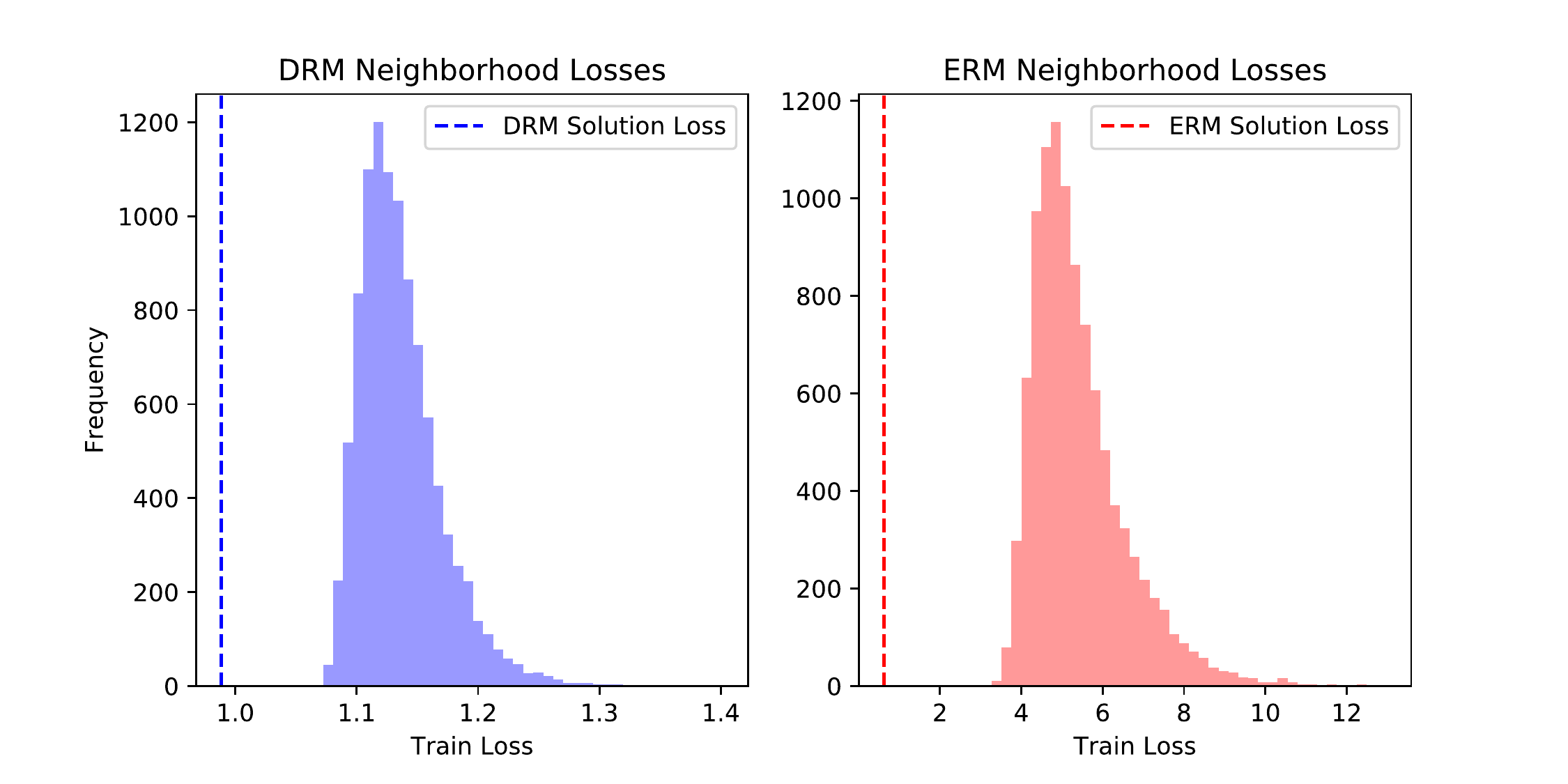}
  \caption{Fully Connected}
    \label{fig:cifar_simple_lipschitz}
    \end{subfigure}
 \begin{subfigure}[b]{.5\linewidth}
  \includegraphics[scale=.4]{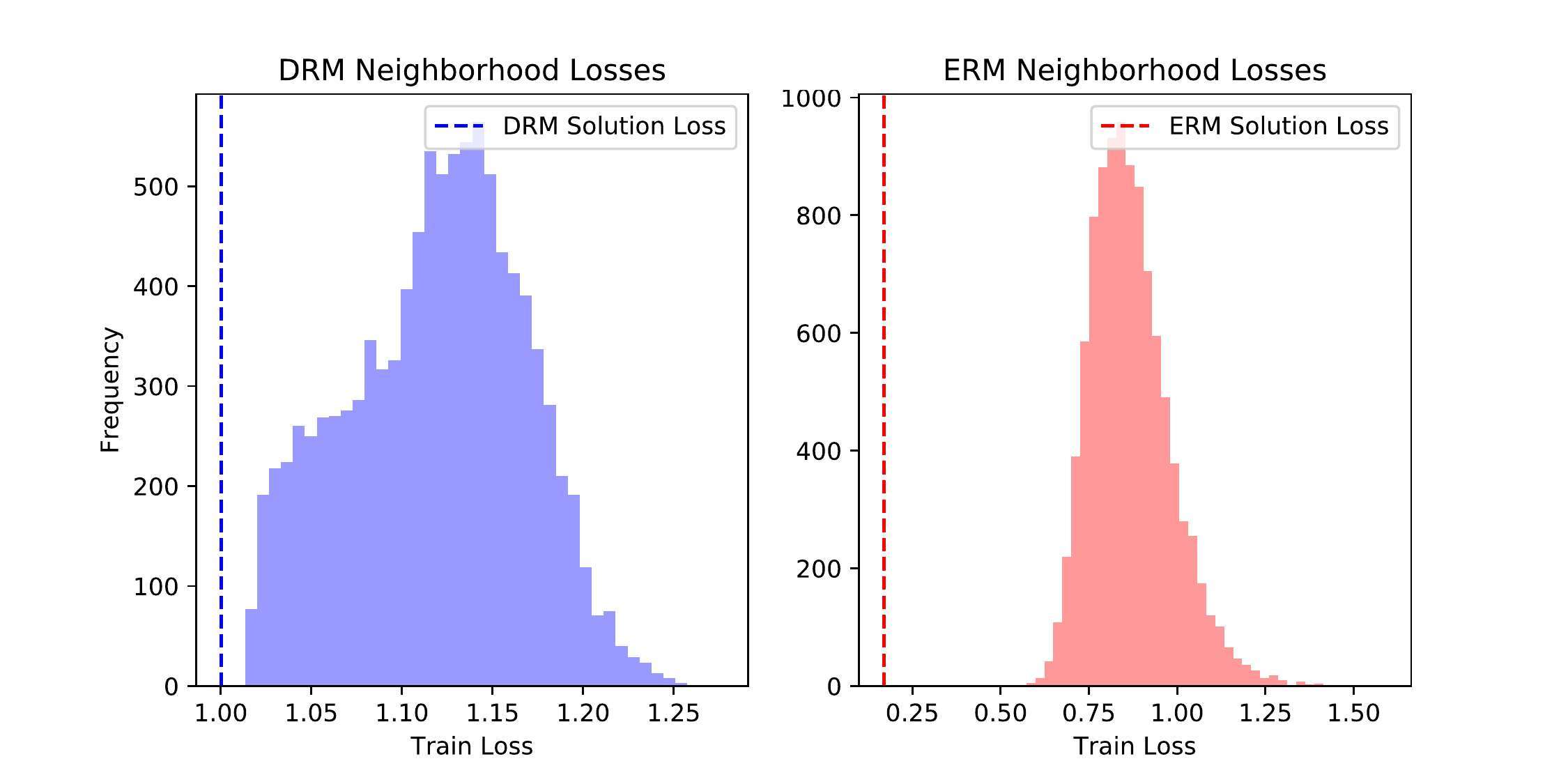}
  \caption{ResNet20}
  \label{fig:cifar_resnet_lipschitz}
      \end{subfigure}
\caption{CIFAR-10 experiments: (a-b) Test accuracy when trained for 800 epochs. (c-d) Distribution of empirical risk for 10,000 random points in $\gamma$-neighborhood of final 800 epoch DRM and ERM solution.}
   \label{fig:cifar}
\end{figure}

\subsection{CIFAR-10}
We perform similar experiments on CIFAR-10 using subclasses airplane, bird, and car and two network architectures: a fully connected and ResNet20. First, we use the same fully connected architecture as before with hidden 3 layers of size (320, 320, 200) with ReLU nonlinearities and an additional fully connected output layer feeding into a 3-class softmax negative log-likelihood objective function. We train the network on a dataset with 50\% of training labels flipped to an incorrect class. For both SGD-ERM and SGD-DRM, we use standard SGD updates with batch size 100 and learning rate .01 for 750 epochs and .001 for the final 50 epochs. For SGD-DRM, the Step 1 perturbations are treated on a layer-wise basis (see Section~\ref{sec:drm_alg_nn}) with $\gamma = 5$. We adopt $r=20$, $q=1$, and a deterministic implementation of Step 6 as above.

Results are nearly identical to those for MNIST; see Figure~\ref{fig:cifar_simple_acc}. SGD-ERM begins by finding a good solution that generalizes, but then continues minimizing the empirical risk and settles into a solution that fits the training data with incorrect labels and thus suffers from a sharp decline in test accuracy. SGD-DRM is again resistant to overfitting and is able to find a good solution that generalizes. Figure~\ref{fig:cifar_simple_lipschitz} illustrates that the DRM solution is in a much flatter portion of the empirical risk landscape than the ERM solution and achieves much smaller diametrical risk at the expense of larger empirical risk. The figure plots the distribution of empirical risk for 10,000 random points in the $\gamma$-neighborhood of the final 800-epoch DRM and ERM solutions. The dashed line represents the value of the empirical risk at the found solution and the rest of the plot is the distribution of empirical risk at points surrounding the solution.

For the ResNet20 architecture, we utilize identical settings except that $\gamma=1$ due to the smaller number of parameters (per layer) and we omit perturbations to the batchnorm layer parameters. Figure~\ref{fig:cifar_resnet_acc} reports test accuracy. Unlike for the fully connect architecture, SGD-DRM suffers nearly zero degradation of test accuracy as training progresses. It also achieves, and maintains, higher test accuracy than is achieved by SGD-ERM. Furthermore, Figure~\ref{fig:cifar_resnet_lipschitz} illustrates again the flatness of the landscape surrounding the SGD-SRM solution. While the diametrical risk is similar for both solutions, there is still a significant gap between the empirical risk of the SGD-ERM solution and points within its neighborhood. For SGD-DRM, this gap is much smaller, indicating a flatter landscape.

\subsection{Discussion}
While these experiments support the proposition that minimization of diametrical risk leads to good generalization and can be used to handle problems with large Lipschitz moduli, there is still much to be explored with DRM. In particular, SGD-DRM can be improved by considering different policies for choosing hyperparameters $\gamma$, $\lambda$, $q$, $d$, batch size, and $p$. The choice of diametrical risk radius $\gamma$, for example, could be chosen adaptively at every iteration along with the step size $\lambda$, mimicking adaptive SGD implementations such as Adam. Also, different choices of $\gamma$ could be made for different groupings of parameters corresponding to the layers in a NN.

An additional algorithmic component that could be made more efficient is the estimation of the diametrical risk, which currently is a bottleneck. First, the sampling in Step 1 can be expensive, especially if done at every iteration or for large value of $d$. Second, if $d$ and/or $q$ is large, the memory requirements of storing $U$ and $V$ can be large. Third, the maximization in Steps 2 and 4 can also be expensive, particularly if the batch size is large and the loss is expensive to compute. Many of these issues, however, can be reduced with parallel implementation. For example, independent workers can each produce a single sample $u_i$ and calculate the value of the objective function. Then, Step 2 can be performed by considering only the collection of $d$ function values produced by the set of workers. We leave these tasks to future work, however, and use the presented results to encourage more work in this direction.\\

\noindent {\bf Acknowledgement.} This work is supported in part by AFOSR under F4FGA08272G001.

\bibliographystyle{abbrv}
\bibliography{refs,refs2}

\end{document}

%% file: mac.tex
\usepackage{graphics}
\usepackage{graphicx}
\usepackage{epsfig}
\usepackage{enumerate}
\usepackage{amsfonts}
\usepackage{theorem}
\usepackage{array}
\usepackage[reqno]{amsmath}
\usepackage{color}
\usepackage{multicol}
\usepackage{subcaption}
\usepackage{float}

\theoremstyle{change}
\newtheorem{proclaim}{PROCLAIM}[section]
\newtheorem{theorem}[proclaim]{Theorem}
\newtheorem{definition}[proclaim]{Definition}

\newtheorem{proposition}[proclaim]{Proposition}

\def\state #1. { \noindent{\bf#1.\enspace}}

\newcommand{\comp}{\,{\raise 1pt \hbox{$\scriptstyle\circ$}}\,}

\newcommand{\exs}{\mathop{\rm exs}}

\newcommand{\prj}{\mathop{\rm prj}}

\newcommand{\reals}{\mathbb{R}}

\newcommand{\natnums}{{{\rm l} \kern -.13em {\rm N} }}
\newcommand{\nats}{\mathbb{N}}
\newcommand{\snats}{{I\kern -.29em N}}
\newcommand{\rats}{{Q\kern -.64em \raise 1pt \hbox{$\scriptstyle |$}\;\,}}
\newcommand{\srats}
    {{Q\kern -.56em \raise 1.2pt \hbox{$\scriptscriptstyle /$}\,}}
\newcommand{\ints}{Z\kern -.46em Z}
\newcommand{\Ex}{\mathbb{E}}
\newcommand{\pluss}{\hskip1pt \raise1pt\vbox{\hrule width6pt \vskip1pt \hrule
                    width6pt} \kern-4pt{\lower1pt\hbox{\vrule height6pt
            \kern1pt\vrule height6pt}}\hskip5pt}
\newcommand{\eop}
    {\hfill{$\vcenter{\hrule height1pt \hbox{\vrule width1pt height5pt
     \kern5pt \vrule width1pt} \hrule height1pt}$} \medskip}
\newcommand{\half}
    {{\raisebox{1pt}{$\frac{1}{2}$}}}

\newcommand{\setd}{{ d \kern -.15em l}}
\newcommand{\hatsetd}{ d \hat{\kern -.15em l }}

\newcommand{\bbP}{\mathbb{P}}

\renewcommand{\epsilon}{\varepsilon}
\renewcommand{\phi}{\varphi}

\newcommand{\tto}{\;{\lower 1pt \hbox{$\rightarrow$}}\kern -12pt
           \hbox{\raise 2.5pt \hbox{$\rightarrow$}}\;}
\newcommand{\overto}[1]{\,{\raise 0pt\hbox{$\rightarrow$}}\kern -9pt
     \hbox{\lower 3pt \hbox{$\scriptscriptstyle#1$}}\hskip6pt}
\newcommand{\underto}[1]{\,{\lower 1pt\hbox{$\rightarrow$}}\kern -9pt
     \hbox{\raise 4pt \hbox{$\,\scriptscriptstyle#1$}}\hskip7pt}
\newcommand{\bigoverto}[1]{{\raise 0pt\hbox{$\,\longrightarrow$}}\kern -16pt
     \hbox{\lower 3pt \hbox{$\scriptscriptstyle#1$}}\hskip4pt}
\newcommand{\bigunderto}[1]{\,{\lower 1pt\hbox{$\longrightarrow$}}\kern -16pt
     \hbox{\raise 4pt \hbox{$\,\scriptscriptstyle#1$}}\hskip6pt}
\newcommand{\bigbigto}[2]{\,{\raise 0pt\hbox{$\,\longrightarrow$}}\kern -16pt
     \hbox{\lower 3pt \hbox{$\scriptscriptstyle#2$}}\kern -10pt
     \hbox{\raise 4pt \hbox{$\,\scriptscriptstyle#1$}}\hskip7pt}
\newcommand{\downto}{{\raise 1pt \hbox{$\scriptstyle \,\searrow\,$}}}
\newcommand{\upto}{{\raise 1pt \hbox{$\scriptstyle \,\nearrow\,$}}}

\newcommand{\notimply}
    {\quad\hbox{$\Longrightarrow \kern -14pt {/}$}\hskip6pt\quad}

\newcommand{\low}[1]{{\lower1pt \hbox{$\scriptstyle #1$}}}
\newcommand{\loww}[1]{{\lower2pt \hbox{$\scriptstyle #1$}}}
\newcommand{\high}[1]{{\raise1pt \hbox{$\scriptstyle #1$}}}

\newcommand{\cA}{{\cal A}}

\newcommand{\nliminf}{\mathop{\rm liminf}\nolimits}

\newcommand{\ninf}{\mathop{\rm inf}\nolimits}
\newcommand{\nsup}{\mathop{\rm sup}\nolimits}

\newcommand{\nnmin}{\mathop{\rm minimize}}
\newcommand{\nargmax}{\mathop{\rm argmax}\nolimits}
\newcommand{\nargmin}{\mathop{\rm argmin}\nolimits}

\newcommand{\lwdy}[2]{\mathrel{\mathop
        {\raisebox{0.1ex}{\null$#1$}}{\hbox{\kern -1.0em
    {\raisebox{-0.8ex}{$\scriptstyle{\;\to #2}$}}}}}}
\newcommand{\lwwdy}[2]{\mathrel{\mathop
        {\raisebox{0.2ex}{\null$#1$}}{\hbox{\kern -1.0em
    {\raisebox{-1.1ex}{$\scriptstyle{\;\to #2}$}}}}}}
\newcommand{\slwwdy}[2]{\scriptsize{{\mathrel{\mathop
        {\raisebox{0.2ex}{\null$#1$}}{\hbox{\kern -1.0em
    {\raisebox{-1.1ex}{$\scriptstyle{\;\to #2}$}}}}}}}}

\def\eto{\,{\lower 1pt\hbox{$\rightarrow$}}\kern -11pt
     \hbox{\raise 4pt \hbox{$\, \scriptstyle e$}}\hskip7pt}
\def\hto{\,{\lower 1pt\hbox{$\rightarrow$}}\kern -11pt
     \hbox{\raise 4pt \hbox{$\, \scriptstyle h$}}\hskip7pt}
\def\pto{\,{\lower 1pt\hbox{$\rightarrow$}}\kern -11pt
     \hbox{\raise 4.5pt \hbox{$\, \scriptstyle p$}}\hskip7pt}
\def\cto{\,{\lower 1pt\hbox{$\rightarrow$}}\kern -11pt
     \hbox{\raise 4pt \hbox{$\, \scriptstyle c$}}\hskip7pt}

\def\dy#1\until#2{\mathrel{\mathop
        {\null#1}\limits^{\hbox{\lower 1.3ex \hbox{$\scriptstyle{\;\to}$}}}}
        {\hbox{{\kern -.2em${\null}
       ^{\hbox{\raise .2ex \hbox{$\scriptstyle{#2}$}}}$}}}}
\def\hidy#1\until#2{\mathrel{\mathop
        {\null#1}\limits^{\hbox{\lower 1.0ex \hbox{$\scriptstyle{\;\to}$}}}}
        {\hbox{{\kern -.2em${\null}
	^{\hbox{\raise 0.9ex \hbox{$\scriptstyle{#2}$}}}$}}}\!}
\def\lody#1\until#2{\mathrel{\mathop
        {\null#1}\limits_{\hbox{\raise 0.5ex \hbox{$\scriptstyle{\to}$}}}}
        {\hbox{{\kern -.25em${\null}
       _{\hbox{\lower 0.7ex \hbox{$\scriptstyle{#2}$}}}$}}}}
\def\lowdy#1\until#2{\mathrel{\mathop
        {\null#1}\limits_{\hbox{\raise 1.0ex \hbox{$\scriptstyle{\to}$}}}}
        {\hbox{{\kern -.25em${\null}
       _{\hbox{\lower 0.9ex \hbox{$\scriptstyle{#2}$}}}$}}}}

\newcommand{\bcdot}{\,{\raise .2ex \hbox{$\cdot$}}\,}
